\documentclass[a4paper, 11pt]{amsart}
\pdfoutput=1
\usepackage{amssymb,amsmath,txfonts,mathrsfs,titletoc,appendix}
\usepackage{graphicx}
\usepackage{latexsym}
\usepackage[alphabetic]{amsrefs}
\usepackage{geometry}
\usepackage{fancyhdr}
\usepackage{color}

 \newtheorem{thm}{Theorem}[section]
 \newtheorem{cor}[thm]{Corollary}
 \newtheorem{lem}[thm]{Lemma}

 \newtheorem{defn}[thm]{Definition}
 \newtheorem{rem}[thm]{Remark}
 \numberwithin{equation}{section}
\newtheorem{lem*}{Lemma}
\newtheorem{cor*}{Corollary} 
\newenvironment{prooff}{\medskip \noindent
{\bf Proof.}}{\hfill \rule{.5em}{1em}
\\}

\newenvironment{pf6}{\medskip \noindent
{\bf Proof of Lemma \ref{sigmaomegacal}}.}{\hfill \rule{.5em}{1em}
\\}

\begin{document}

\title{On extending calibration pairs}

\author{Yongsheng Zhang}
\address{\em{Currently visiting:} 
Max Planck Institute for Mathematics,
Vivatsgasse 7,
53111 Bonn,
Germany.}
\email{yongsheng.chang@gmail.com}
\date{\today}
\date{\today}
\thanks{This work was partially sponsored by
NSFC (Grant Nos. {11526048, 11601071}),
the Fundamental Research Funds for the Central Universities, the SRF for ROCS, SEM
and 
the NSF under Grant No. 0932078 000,
while the author was in residence at the MSRI during the 2013 Fall.}
\begin{abstract}
The paper studies how to extend local calibration pairs to global ones in various situations.
As a result, new discoveries involving mass-minimizing properties are exhibited.
In particular, we show that a $\mathbb R$-homologically nontrivial connected submanifold $M$
of a smooth Riemannian manifold $X$
is homologically mass-minimizing for some metrics in the same conformal class.
Moreover, several generalizations for $M$ with multiple connected components or for a mutually disjoint collection
(see \S \ref{sr}) are obtained.
For a submanifold with certain singularities, we also establish an extension theorem for generating global calibration pairs. 
By combining these results,
we find that, in some Riemannian manifolds,
there are homologically mass-minimizing smooth submanifolds
which cannot be calibrated by any smooth calibration.
\end{abstract}
\maketitle
\titlecontents{section}[0em]{}{\hspace{.5em}}{}{\titlerule*[1pc]{.}\contentspage}
\titlecontents{subsection}[1.5em]{}{\hspace{.5em}}{}{\titlerule*[1pc]{.}\contentspage}
{\setcounter{tocdepth}{2} \small \tableofcontents}
\section{Introduction}\label{Section1} 

             A calibration on a smooth Riemannian manifold $(X,g)$ is a closed differential $m$-form $\phi$
          whose value at every point on every unit $m$-plane is at most one.
              The fundamental theorem of calibrated geometry \cite{HL2} asserts that 
          an $m$-dimensional oriented compact submanifold $M$ (or more generally a current) of $X$
          for which $\phi$ has value one a.e. on every unit tangent plane is mass-minimizing in its homology class of normal currents. 
         We call $(\phi, g)$ a calibration pair of $M$ on $X$.
         
         In this paper we shall create such balanced pairs for objects in various situations.
         The idea is to have a local calibration pair and extend it to a global one.
         Based on types of objects to deal with, the paper divides into two parts: the smooth case and the singular case.
         
         Given a homologically nontrivial, oriented, connected, compact submanifold $M$,
         we show that one can conformally change the metric on $X$ such that
         $M$ becomes homologically mass-minimizing.
         Our existence result in every conformal class of metrics generalizes the existence theorem
         of Tasaki \cite{Tasaki}.
         In his paper, Tasaki first applied a functional analysis argument of Sullivan \cite{S}
         for a global form $\phi$ which has positive values on the oriented tangent planes of $M$,
         and then he made use of two powerful results (Lemmas \ref{hl1} and \ref{hl2} in our paper)
         by Harvey and Lawson \cite{HL1} to build a metric $g$ adapted to $\phi$
         so that $(\phi, g)$ forms a calibration pair of $M$.
            Although our improvement reflects in the metric slot,
         the essential difference comes from the construction of calibrations.
         Our recipe is this.
         Lemma \ref{fperp} provides a well-behaved local calibration pair.
         We first extend the form to a global one in \S\ref{forms},
         and then glue metrics accordingly in \S\ref{metrics} for a global calibration pair. 
         The case of a constellation of mutually disjoint submanifolds possibly of different dimensions is also studied in \S\ref{sr}.

        By Federer and Fleming \cite{FF}
        there exists at least one mass-minimizing normal current in every real-valued homology class of a compact Riemannian manifold.
        However the regularity of these mass-minimizing currents and their distributions are quite complicated in general. 
        We construct nice metrics in Theorem \ref{2} so that (as functionals over smooth forms) 
        all homologically mass-minimizing currents of codimension at least 3 are just linear combinations
        of (integrations over) submanifolds.
        The idea is the following. 
        For each dimension, the homology space has a basis that can be represented by oriented connected compact submanifolds.
        One can arrange these representatives so that all intersections among them are transversal.
        Then enough calibrations can be constructed to feed our need for codimensions no less than 3.
        
        Except in low dimensions, mass-minimizing and even calibrated currents of codimension one can have singularities.
        N. Smale \cite{NS} gave the first examples of homologically mass-minimizing compact hypersurfaces with isolated singular points.
        In the second part of our paper, 
        a different method for getting such creatures through calibrations is gained.
        We first establish an extension result
        Theorem \ref{conecal} (also see Example 1)
        which allows us to extend a ``nice" local calibration pair
\footnote{In this case the calibration form may be singular somewhere. See the definition of coflat calibration in \S\ref{P}.}
         of a singular submanifold around its singular set to a calibration pair on some neighborhood of the submanifold.
        Under certain conditions a further extension to a global pair can be made.
       Then in Example 2 we illustrate how to build up examples satisfying the requirements in Theorem \ref{conecal}.
       They provide lots of examples similar to N. Smale's.
        
        Our local models of singularities with nice calibration pairs
        include all homogeneous mass-minimizing hypercones
         which have (coflat) calibrations singular only at the origin (see \cite{Z4}),
        and all special Lagrangian cones 
        (see
\cite{Joyce, Mc, CM, Haskins, HK, HK2, HK3})
          that enjoy smooth calibrations.
        {In fact, 
        we show in \cite{Z3} that
        every
          {known}
          area-minimizing hypercone and every oriented area-minimizing cone in \cite{Law, TZ, XYZ2}
        can be realized as a tangent cone at a singular point of some homologically area-minimizing singular compact submanifold.}
          
       A very interesting phenomenon, that we observe in Example 3, 
        is the existence of homologically mass-minimizing smooth submanifolds
        which cannot be calibrated by any smooth calibration.
        Actually, all coflat calibrations of the submanifold share at least one common singular point.
        By Remark \ref{88} there are examples for which calibrations share more complicated singular sets.
        
        Through blowing-up we get Example 4 which relates to twisted calibrations  \cite{TM} and integral currents mod 2 \cite{Ziemer}.
        It gives us a non-orientable compact singular hypersurface that is mass-minimizing in its homology class of integral currents mod 2.
        
        {\ }
        
$\text{\sc Acknowledgement}.$
This paper is an expansion of part of the author's Ph.D. thesis at Stony Brook University.
He wants to express particular thanks to 
Professor H. Blaine Lawson, Jr. for his guidance and constant encouragement.
He also wishes to thank Professor Frank Morgan for several valuable communications,
Professor Xiaobo Liu for the informative conversation,
Professor Brian White for his help on Example 4 of \S \ref{exs},
and the MSRI for the warm hospitality.

{\ }


\section{Preliminaries}\label{P} 
            We review some fundamental concepts and results in calibrated geometry.
            Readers are referred to \cite{HL2} for a further understanding on this subject
            and to \cite{FM} for a quick overview of geometric measure theory.
            
            \begin{defn}\label{comassdef}
            Let $\phi$ be a smooth $m$-form on a Riemannian manifold $(X,g)$.
            At a point $x\in X$ we define the \textbf{comass} of $\phi_x$ to be
                     \begin{equation*}
                     \|\phi\|_{x,g}^*=\max \ \{\phi_x( \overrightarrow V_x) : \overrightarrow V_x \ \text{is a unit simple m-vector at x}\}.
                     \end{equation*}
             Here {``simple"} means $\overrightarrow V_x=e_1\wedge e_2\cdots \wedge e_m$
             for some $e_i\in T_xX$.
             \end{defn}

            \begin{rem}\label{comassf}
            $\|\phi\|_{g}^*$ will be viewed as a pointwise function in this paper.
            In general it is merely continuous. 
            At a point $x$ where $\phi_x\neq0$,
                       \begin{equation*}
                            \begin{split}
                            \|\phi\|^*_{x,g}& =  \max \{\phi(\overrightarrow V_x): \overrightarrow V_x\ \text{is a simple $m$-vector  at } x \text{ with}\ 
\|\overrightarrow V_x\|_{g}=1\}\\
                  &=\max\{1/{\|\overrightarrow V_x\|_{g}}: \overrightarrow V_x\ \text{is\ a\ simple\ $m$-vector\ at } x \text{ with}\ \phi(\overrightarrow V_x)=1\}\ \\
                  &=1/\min\{\|\overrightarrow V_x\|_{g}:\overrightarrow V_x\ \text{is\ a\ simple\ $m$-vector\ at } x \text{ with}\ \phi(\overrightarrow V_x)=1\}.
                             \end{split}
                         \end{equation*}
             \end{rem}
             
               \begin{defn}
               Denote the dual complex of the {\em de Rham} complex of $X$ by $(\mathscr E'_*(X),d)$.
               Elements of $\mathscr E'_k(X)$ are $k$-dimensional {\em de Rham} \textbf{currents} (with compact support)
               and $d$ is the adjoint of exterior differentiation.
               \end{defn}
             
               \begin{defn}\label{dmass}
               In $(X,g)$, the \textbf{mass} $\mathrm{\mathbf{M}}(T)$ of $T\in\mathscr E'_k(X)$ is defined to be
               $$\sup \{T(\psi):\psi\ smooth\ m\text{-}form\ with\  \sup_X\|\psi\|_{g}^*\leq 1\}.$$
               \end{defn}
          When $\mathrm{\mathbf{M}}(T)<\infty$, $T$ determines a unique {Radon} measure $\|T\|$ characterized by
                 $$ \int_X f\cdot d\|T\|=\sup\{T(\psi): \|\psi\|_{x,g}^*\leq f(x) \}$$
          for any nonnegative continuous function $f$ on $X$.
          Therefore $\mathrm{\mathbf{M}}(T)=\|T\|(X)$.
          Moreover, the {Radon}-{Nikodym} Theorem asserts the existence of
          a $\|T\|$ measurable tangent $m$-vector field $\overrightarrow T$ a.e. with 
          vectors $\overrightarrow T_x \in \Lambda^m T_xX$ of unit length in the dual norm of the comass norm, satisfying
                \begin{equation}\label{current}
                T(\psi)= \int_X\psi_x(\overrightarrow {T_x})\ d \|T\|(x)\ \  \text{for any smooth $m$-form }\psi,
                \end{equation}
           \text{or\ briefly}
           $T = \overrightarrow T\cdot \|T\|\ a.e.\ \|T\|.$
           When $T$ has local finite mass, one can get {Radon} measure $\|T\|$ and decomposition (\ref{current}) as well.
                 \begin{defn}
                 For a function $f$,
                 set $\mathbf{spt}(f)$ to be its support.
                 For a current $T$, let $U_T$ stand for the largest open set with $\|T\|(U_T)=0$.
                 Then the support of $T$ is denoted by $\mathbf{spt}(T)=U_T^c$.
                 \end{defn}
                 \begin{defn}
                 Let $\mathbb M_k(X)=\{T\in\mathscr E'_k(X): \mathrm{\mathbf{M}}(T)<\infty\}$.
                 Then $N_k(X)=\{T\in\mathbb M_k(X): dT\in\mathbb M_{k-1}(X)\}$ is the space of $k$-dimensional \textbf{normal\ currents}.
                 \end{defn}
                 \begin{rem}
                 We view a current in $\mathbb M_k$ as a functional over smooth $k$-form not a specific representative of generalized distribution.
                 \end{rem}
                 Note that $(N_*(X),d)$ form a chain complex.
                 Recalling the natural isomorphisms established by {de Rham}, Federer and Fleming:
                            $$H_*(\mathscr E'_*(X))\cong H_*(X;\mathbb R)\cong H_*(N_*(X))$$
                 we identify these three homology groups.

           \begin{defn}\label{calibration}
           A smooth form $\phi$ on $(X,g)$ is called a \textbf{calibration} if 
           $\sup_{X}\|\phi\|_{g}^*= 1$ 
           and
           $d\phi=0.$
           Such a triple $(X,\phi,g)$ is called a \textbf{calibrated manifold}. 
                   If $M$ is an oriented submanifold with $\phi|_{M}$ equal to the volume form of $M$,
                   then $(\phi,g)$ is a \textbf{calibrated pair} of $M$ on $X$.
          We say $\phi$ \textbf{calibrates} $M$ and $M$ \textbf{can be calibrated} in $(X,g)$.
          \end{defn}
          
          \begin{defn}\label{calibratable}
                   Let $\phi$ be a calibration on $(X,g)$.
                   We say that a current $T$ of local finite mass is \textbf{calibrated} by $\phi$, if 
                        $\phi_x(\overrightarrow T_x)=1\ a.a.\ x\in X\ \text{for}\ \|T\|.$
          \end{defn}

                  \begin{rem}
                               For an oriented compact submanifold $M$,
                               the current $[[M]]=\int_M\cdot\ $ is calibrated if and only if
                               $M$ is calibrated.
                               \end{rem}

                    The following is the fundamental theorem of calibrated geometry in \cite{HL2}.
                              \begin{thm}\label{hl}
                                           If $T$ is a calibrated current 
                                \footnote{It is called a positive $\phi$-current in \cite{HL2}.}
                                with compact support
                                           in $(X,\phi,g)$ and
                                           $T'$ is any compactly supported current homologous to $T$(i.e., $T-T'$ is a boundary and in particular $dT=dT'$),
                                           then
                                                          \begin{equation*}
                                                          \mathrm{\mathbf{M}}(T)\leq  \mathrm{\mathbf{M}}(T')
                                                          \end{equation*}
                                          with equality if and only if $T'$ is calibrated as well.
                               \end{thm}
                               
                               It is often useful to allow calibrations to have certain singularities.
                               
                                \begin{defn}\label{coflat}
           Let $\phi$ be a calibration of degree $m$ on $X-S_\phi$,
           where $S_\phi$ is a closed subset of $X$ of Hausdorff $m$-measure zero.
           Then $\phi$ is called a \textbf{coflat calibration} on $X$.
           We say $\phi$ \textbf{calibrates} a current, if it is calibrated by $\phi$ on $X-S_\phi$.
          \end{defn}

              \begin{rem}\label{coflatFTCG}     
               Actually there is a coflat version (Theorem 4.9 in \cite{HL2}) of the fundamental theorem of calibrated geometry,         
                     and a current calibrated by a coflat calibration is homologically mass-minimizing as well.
                \end{rem}

 \section{Smooth case}\label{SC}               
        We shall use some properties of comass. 
        Especially, {Lemma \ref{fperp}} is crucial to our method and {Lemma \ref{CCGP}} provides certain control on comass while gluing metrics.

              \subsection{Properties of comass}\label{poc}
                        \begin{lem}\label{norm}
                              For any metric $g$, $m$-form $\phi$ and positive function $f$ on $X$,
                                   $$\|\phi\|_{f\cdot g}^*=f^{-\frac{m}{2}}\cdot \|\phi\|_g^*.$$
                        \end{lem}
                             \begin{prooff} By the formula in {Remark \ref{comassf}}.\end{prooff}
                        \begin{lem}\label{big}
                        For any $m$-form $\phi$ and metrics $g'\geq g$ on $X$, we have
                        $$\|\phi\|_{g'}^*\leq\|\phi\|_g^*.$$ 
                        \end{lem}
                               \begin{prooff}
                               By the definition of comass.
                               \end{prooff}
                         \begin{lem}[Comass control for gluing procedure]\label{CCGP}
                         For any $m$-form $\phi$, positive functions $a$ and $b$, and metrics $g_1$ and $ g_2$,
                         it follows
                                 \begin{equation}\label{gluemetrics}
                               \|\phi\|^*_{ag_1+bg_2}\leq\frac{1}{\sqrt{a^m\cdot \frac{1}{\|\phi\|^{*2}_{g_1}}+b^m\cdot \frac{1}{\|\phi\|^{*2}_{g_2}}}}
                                 \end{equation}
                         where $\frac{1}{0}$ and $\frac{1}{+\infty}$ are identified with $+\infty$ and $0$ respectively.
                         \end{lem}
                          \begin{prooff}
                          The statement is trivial at where $\phi$ vanishes.
                          Suppse $\phi_x\neq0$ at a point $x$.
                          In the subspace spanned by a simple $m$-vector $\overrightarrow V_x$,
                          there exists an orthonormal basis $(e_1,\ \cdots, e_m)$ of $g_1$,
                          under which $g_2$ is diagonalized as $diag(\lambda_1,\ \cdots ,\lambda_m)$ for some $\lambda_i>0$.
                          Let $\overrightarrow V_x=te_1\wedge\cdots\wedge e_m$, then
                                     \begin{equation}\label{inq}
                                            \begin{split}
                                              \|\overrightarrow V_x\|_{ag_1+bg_2}^2 &= t^2(a+b\lambda_1)\cdots(a+b\lambda_m)\\
                                                                                                            &=\  \ t^2[a^m+\cdots+b^m{\bold \Pi}\lambda_i]\\
                                                                                                            &\geq \ \ \quad t^2a^m+t^2b^m{\bold\Pi}\lambda_i\\
                                                                                                            &= a^m\|\overrightarrow V_x\|_{g_1}^2+b^m\|\overrightarrow V_x\|_{g_2}^2. 
                                              \end{split}
                                       \end{equation}
                            By {Remark \ref{comassf}},
                           (\ref{inq}) implies \eqref{gluemetrics}.
                           \end{prooff}
                          \begin{lem}[Comass one lemma]\label{fperp}
                          Suppose $(E, \pi)$ is a disk bundle over an oriented manifold $M$ (as the zero section) and $g$ is a Riemannian metric on $E$.
                          Then each fiber is perpendicular to $M$ if and only if $\pi^*\omega$ has comass one pointwise along $M$
                          where $\omega$ is the induced volume form of $M$. 
                          \end{lem}
                        \begin{prooff}
                        For $x\in M$, take an oriented orthonormal basis $\{e_1,\ \cdots,\ e_m\}$ of $T_xM$. 
                        Then we have unique decompositions
                        $e_i=\sin\theta_i\cdot a_i+ \cos\theta_i\cdot b_i$
                        where $b_i$ is some unit vector in $F_x$ $-$ the subspace of fiber directions in $T_xE$,
                        $a_i$ is a unit vector perpendicular to $F_x$ ,  
                        and $\theta_i$ is the angle between $e_i$ and $F_x$. 
                        By the choice of $\{e_i\}$,
                               \begin{equation}\label{q1}
                                    \begin{split}
                                         {\ } \ \ \ \ \ \ \ \ 1&=\omega(e_1\wedge e_2\cdots\wedge e_m)\\
                                                                &=\pi^*\omega(e_1\wedge e_2\cdots\wedge e_m)\\
                                                                &=\pi^*\omega(\sin\theta_1\cdot a_1\wedge\cdots\wedge\sin\theta_m\cdot a_m)\\
                                                                &={\bold \Pi} \sin\theta_i\cdot \pi^*\omega(a_1\wedge a_2\cdots\wedge a_m).
                                     \end{split}
                                \end{equation}
                        The third equality is because that elements of $F_x$ annihilate $\pi^*\omega$.
                        Since $\{a_i\}$ are of unit length,
                              $\|\pi^*\omega\|^*_{x,g}\geq 1, \ \forall x\in M$,
                        with equality if and only if $F_x \perp T_xM$.
                       \end{prooff}
\begin{rem}\label{loc}
                      Since $\pi^*\omega$ is smooth and simple, $\|\pi^*\omega\|^*_g$ is smooth.
                      When fibers are perpendicular to $M$,
                      $(\pi^*\omega, (\|\pi^*\omega\|^*_g)^{\frac{2}{m}}g )$
                      gives a calibration pair of $M$ on $E$.
  \end{rem}
{\ }
\subsection{Global forms}\label{basic}
                In the singular homology theory the {\em Kronecker} product $<\cdot ,\cdot >$ between cochains and chains induces a homomorphism
                        \begin{equation*}\label{kappa}
                        \kappa\ : H^q(X; G)\rightarrow \ \text{Hom}_{\mathbb Z}(H_q(X;\mathbb{Z}),\ G) \text{ given by}
                        \end{equation*}
                        \begin{equation*}
                        \kappa\ ([z^q])([z_q])\triangleq <[z^q],[z_q]>
                        \end{equation*}
                        where $G$ is an abelian group.
                        A classical result asserts that
                       $\kappa$ is surjective. 
                       When $G=\mathbb{R}$, by the {de\ Rham} {Theorem},
                       $\kappa\ : H_{dR}^q(X)\twoheadrightarrow \  \text{Hom}_{\mathbb{R}}(H_q(X;\mathbb{R}),\mathbb{R}).$
                       
                Suppose $\{M_\alpha\}$ are mutually disjoint $m$-dimensional oriented connected compact submanifolds 
                with homology classes $\{[M_\alpha]\}$ lying in one common side of some hyperplane through the zero of $H_m(X;\mathbb{R})$.
                Then there exists a homomorphism 
                $\digamma\in \text{Hom}_{\mathbb{R}}(H_m(X;\mathbb{R}),\mathbb{R})$
                forwarding $\{[M_\alpha]\}$ to positive numbers.
                As a consequence, we have the following.
                
                        \begin{lem}\label{formmore}
                        Suppose 
                        $\{M_\alpha\}$ satisfy the above condition.
                        Then there exists a closed $m$-form $\phi$ on $X$ with
                        $\int_{M_\alpha} \phi>0\ \text{ for  each } M_\alpha.$
                        \end{lem}
                        

{\ }
\subsection{Gluing of forms}\label{forms}
                       Given an oriented connected compact submanifold $M$ in $(X,g)$,
                       consider its $\epsilon$-neighborhood $U_\epsilon$.
                       When $\epsilon$ is small enough, 
                       the metric induces a disk bundle structure of $U_\epsilon$,
                       whose fiber is given by the exponential map restricted to normal directions of $M$.
                       Hence by Remark \ref{loc} a local calibration pair of $M$ can be produced.
                       We shall extend (a modification of) this local pair to a global one.
                       Let us glue forms first.

                       By a strong deformation retraction from $U_{\epsilon}$ to $M$,
                       $H^m(U_\epsilon;\mathbb{R})\cong H^m(M;\mathbb{R}).$
                       Therefore for any $[\phi_1]$ and $ [\phi_2]\in$ $H^m(U_\epsilon;\mathbb{R})$
                               \begin{equation}\label{equt}
                                     [\phi_1]=[\phi_2]\ \Leftrightarrow \int_M \phi_1=\int_M \phi_2.
                               \end{equation}
                               
                    Assume further $[M]\neq[0]\in H_m(X;\mathbb R)$.
                    By {\S \ref{basic}} there exists a closed $m$-form $\phi$ on $X$ with $s=\int_M \phi> 0.$
                    Let $\pi$ be the projection map of the disk bundle.
                    Then in $U_\epsilon$ 
                           \begin{equation}\label{ef} 
                           \int_M \frac{s\cdot\pi^*\omega}{\text{Vol}_g(M)}=s=\int_M  \phi
                           \end{equation}
                    where $\omega$ is the volume form of $M$.
                    Denote the integrand of the left hand side of \eqref{ef} by $\omega^*$.
                    By (\ref{equt}) $[\omega^*]=[\phi]$ in $H^m(U_\epsilon;\mathbb{R})$ which indicates
                              \begin{equation}\label{abouttoglue}
                              \phi=\omega^*+d\psi
                              \end{equation}
                     for some smooth $(m-1)$-form $\psi$ on $U_\epsilon$.
                     Now take $\Phi=\omega^*+d((1-\rho(\textbf{d}))\psi)$
                     where \textbf{d} is the distance function to $M$ and $\rho$ is given in the picture.
 \begin{figure}[ht]
\centering
\includegraphics[scale=0.21]{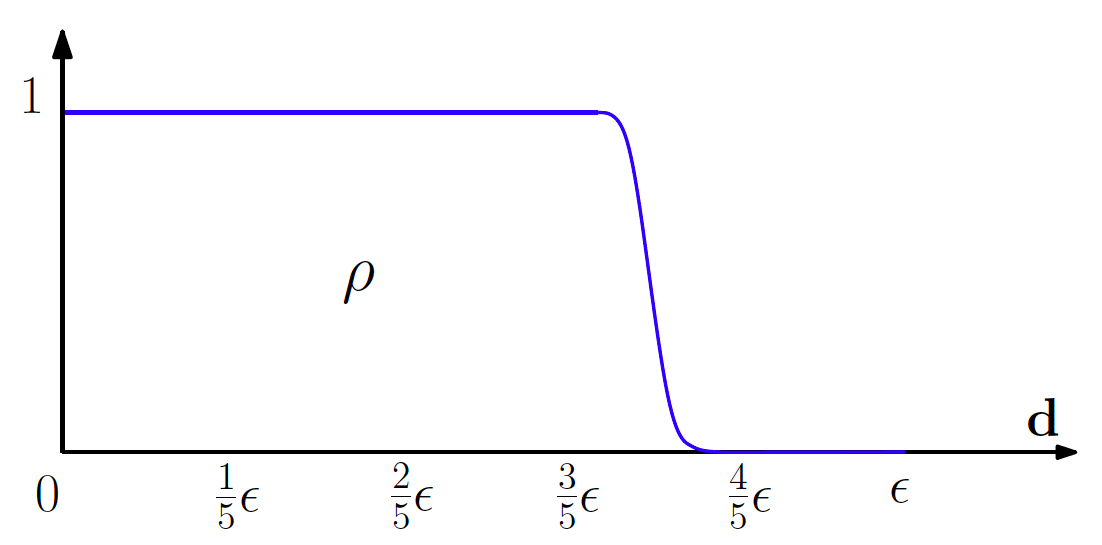}
\end{figure}
                      Clearly $\Phi$ extends to a closed smooth form on $X$:
                      \begin{equation*}
                                \Phi=\begin{cases}\label{glueforms}
                                          \omega^* & 0\leq  \textbf{d}\leq \frac{3}{5}\epsilon\\
                                          \omega^* + d((1-\rho( \textbf{d}))\psi)\ \ \ \ \ \ \ & \frac{3}{5}\epsilon <  \textbf{d} \leq \frac{4}{5}\epsilon\\
                                          \phi & \frac{4}{5}\epsilon< \textbf{d}\\
                                           \end{cases}
                       \end{equation*}
                      By Remark \ref{loc} and Lemma \ref{poc}, $(\Phi, (\|\Phi\|^*_g)^{\frac{2}{m}}g )$ is a calibration pair of $M$ for $\textbf{d}< \frac{3}{5}\epsilon$.

{\ }

\subsection{Gluing of metrics}\label{metrics}
                          Our goal is to extend $(\|\Phi\|^*_g)^{\frac{2}{m}}g $ to a global metric
                          under which the global form $\Phi$ becomes a calibration. 
                          Choose an appropriate positive smooth function $\alpha$ such that
                                 \begin{equation}\label{globalcontrolofalpha}
                                 \|\Phi\|^*_{\alpha g}<1\text{ on } X,
                                 \end{equation}
                            and a gluing function $\sigma=\sigma$(\textbf{d}) shown in the picture.
              \begin{figure}[h]
              \begin{center}
              \includegraphics[scale=0.26]{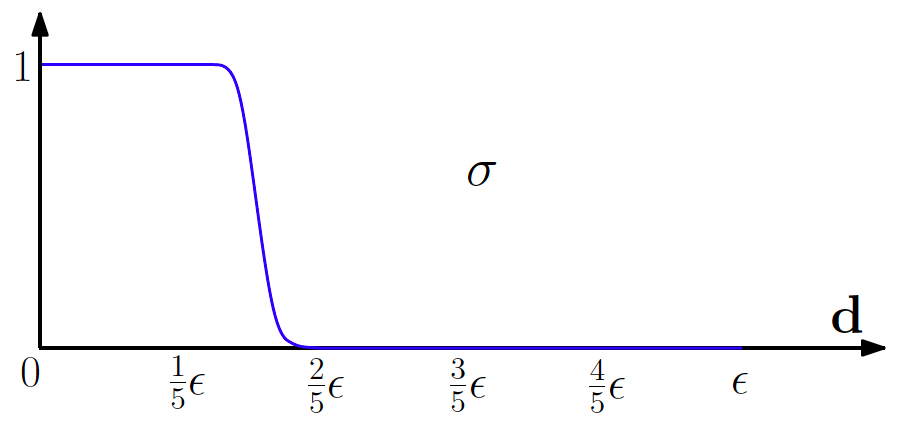}
              \end{center}
              \end{figure}
              Then by Lemmas \ref{poc}-\ref{CCGP}
              \begin{equation}\label{mgofcc}
                      \tilde g=\sigma^{\frac{1}{m}}(1+ \textbf{d}^2)(\|\Phi\|^*_g)^{\frac{2}{m}}g +\alpha(1-\sigma)^{\frac{1}{m}}g
              \end{equation}
              can serve for our purpose.
              Here the factor $(1+ \textbf{d}^2)$ makes $\bold{spt}(\|\Phi\|^*_{\tilde g}-1)=M$
              which implies that $[[M]]$ is uniquely mass-minimizing in $[M]$.

{\ }

\subsection{Some results}\label{sr}

             We can have a few consequences of the constructions in \S\ref{forms} and \S\ref{metrics}.
             An immediate one is this.

             \begin{thm}\label{mainthm2}
             Suppose $(X,g)$ is a Riemannian manifold and
             $M$ is an oriented connected compact $m$-dimensional submanifold
             with $[M]\neq0\in H_m(X;\mathbb{R})$.
             Then there exists a metric $\hat g$ conformal to $g$
             such that $[[M]]$ is the unique mass-minimizer in its homology class in $(X, \hat g)$.
             \end{thm}
            \begin{rem}\label{cptcst}
            When $X$ is compact, $\alpha$ in \eqref{globalcontrolofalpha} can be a sufficiently large constant.
            Set $\hat g\triangleq  \alpha^{-1}\tilde g$ and $\hat \Phi \triangleq  \alpha ^{-\frac{m}{2}}\Phi$.
            Then, by Lemma \ref{poc}, $(\hat \Phi, \hat g)$ is a calibration pair of $M$ and $\hat g=g$ on $X-U_\epsilon$.
            \end{rem}
            Assume that $M$ is an oriented submanifold with (countably many) connected components $\{M_i\}$
            and that every $M_i$ is compact.
            If $\{M_i\}$ satisfies the condition in \S\ref{basic},
            then the same procedure works and we have the following.
             \begin{thm}\label{many}
             Let $M$ be given as above in $(X,g)$.
             Then there exist a metric $\hat g$ conformal to $g$ and a calibration $\hat \Phi$
             such that 
             every nonzero current $T= \sum_{i}t_i[[M_i]]$,
             where $\{t_i\}$ are nonnegtive and only finitely many of them are nonzero, is calibrated in $(X,\hat \Phi, \hat g)$.
             \end{thm}
          \begin{rem}
             Suppose an oriented compact submanifold $M$ has exactly two connected components $M_1$ and  $M_2$ with $[M_1]=[M_2]$.
             Then under $\hat g$ all mass-minimizing currents in the class are given by 
               $t[[M_1]]+(1-t)[[M_2]]$ where $0\leq t\leq 1$.
             \end{rem}            
            When each $[M_i]$ is nonzero,
            one can choose some hyperplane $\mathcal{P}_m$ through zero in $H_m(X;\mathbb{R})$
            that avoids all classes $\{[M_i]\}$. 
            Now $\mathcal{P}_m$ divides the space into two open chambers.
            By reversing orientations of components in one chamber,
            we get a new collection satisfying the requirement in \S\ref{basic}.
           \begin{cor}
             Suppose each $[M_i]$ is nonzero.
             Then in any conformal class of metrics there exists a metric $\hat g$
             such that 
             every $[[M_i]]$ is homologically mass-minimizing in $(X, \hat g)$.
             \end{cor}        
              In order to have a clearer description in more general situation, we introduce some notations.
              \begin{defn}\label{collection}
              A family $\mathfrak M$ of mutually {disjoint} oriented connected compact submanifolds of $X$
              is called a \textbf{mutually disjoint collection}
              and an element of $\mathfrak M$ is a \textbf{component}.
              The (nonempty) subset $\mathfrak M_k$ of all components of dimension $k$ is its \textbf{k\text-level}.
              \end{defn}
              \begin{defn}\label{neat}
              Let $\mathfrak M=\{M_i\}_{i=1,2,\cdots}$ be a mutually disjoint collection of countably many components.
              If the set $\bigcup_{i\neq j}M_i$ is closed for every $j$, 
              then $\mathfrak M$ is called a \textbf{neat collection}.
              \end{defn}
              
              The neatness implies the existence of $\epsilon_i>0$ such that $\{U_{\epsilon_i}(M_i)\}$ are mutually disjoint.
              
              \begin{thm}\label{Xnoncthm}
              Suppose that $\mathfrak M$ is a neat collection
              and that each component represents a nonzero class in the $\mathbb{R}$-homology of $X$.
              In addition, assume every level except the lowest of $\mathfrak M$ has finite components. 
              Then in any conformal class of metrics there exists a metric $\hat g$
             such that 
             each $[[M_i]]$ is homologically mass-minimizing in $(X, \hat g)$.
             \end{thm}

                      \begin{prooff}
                            Without loss of generality, let $\mathfrak M=\{A^a,B^b\}$ with $a>b$ and $g$ be a metric. 
                            Take small positive $\epsilon_{1,2}$ for the procedure in \S\ref{forms}
                            so that $U_{\epsilon_1}(A)$ and $U_{\epsilon_2}(B)$ are disjoint. 
                            Suppose one gets an $a$-form $\Phi$ for $A$.
                            Then $\Phi=d\theta$ in $U_{\epsilon_2}(B)$ for some form $\theta$ of degree $a-1$. 
                            So $\Phi$ can be assumed identically zero in $U_{\epsilon_2}(B)$ from the beginning.
                            Using an $\alpha_a$ whose value remains one on $U_{\epsilon_2}(B)$
                            we get a metric $\tilde g$ by {\S \ref{metrics}} under which $A$ is calibrated by $\Phi$.
                            
                            By the compactness of $\overline{U_{\epsilon_1}(A)}$,
                            there is a $b$-form $\psi$ with $\int_B\psi>0$ and $\|\psi\|_{\tilde g}^*<1$ on $\overline{U_{\epsilon_1}(A)}$.
                            Suppose we get $\Psi$ following \S\ref{forms}.
                            Then one can find an $\alpha_b\geq 1$ 
                            with value one in $\overline{U_{\epsilon_1}(A)}$ for           
                                 $$ \|\Psi\|^*_{\alpha_b \tilde g}<1\text{ on } X.$$
                            By \S\ref{metrics} we get a calibration pair $(\Psi,\hat g)$ of $B$.
                            Note that $(\Phi, \hat g)$ is a calibration of $A$.
                            \end{prooff}
                            \begin{rem}
                            The compactness of $\overline{U_{\epsilon_1}(A)}$ is important.
                            If a level of $\mathfrak M$ has infinitely many components, then
                            our current proof cannot descend further from that level.
                           \end{rem}

                          In \cite{Tasaki} Tasaki studied the ``equivariant" case.
                          \begin{thm}[Tasaki] 
                          Let $K$ be a connected compact Lie transformation group of a manifold $X$ and
                          $M$ be a (connected) compact oriented submanifold in $X$.
                          Assume $M$ is invariant under the action of $K$
                          and it represents a nonzero $\mathbb R$-homology class of $X$.
                          Then there exists a $K$-invariant Riemannian metric $g$ on $X$ such that
                          $M$ is mass-minimizing in homology class with respect to $g$.
                          \end{thm}
                          
                          By our method, one can improve the result.
                          \begin{thm}\label{eqcase2}
                          Let $K$ be a compact Lie transformation group of a manifold $X$ and
                          $M$ be a connected compact oriented submanifold with $[M]\neq 0\in H_m(X;\mathbb R)$.
                          Assume $M$ is invariant under the action of $K$ and the action is orientation preserving.
                          Then for any $K$-invariant Riemannian metric $g^K$,
                          there exists a $K$-invariant metric $\hat g^K$ conformal to $g^K$
                          such that $M$ can be calibrated in $(X,\hat g^K)$.
                          \end{thm}
                          
                          \begin{prooff}
                          There is a {Haar}-measure $d\mu$ with $\int_Kd\mu=1$ for compact $K$.
                           Since the action is orientation preserving and $g^K$ is $K$-invariant,
                           $\omega^*$ and $\textbf{d}$ are $K$-invariant.
                           So one can use $d\mu$ to average \eqref{abouttoglue} for a $K$-invariant $\Phi$
                           which equals $\omega^*$ in $M$.
                            Then average the corresponding $\alpha$.
                            By \eqref{mgofcc} one can get a $K$-invariant calibration pair
                            $(\Phi, \hat g^K)$.
                            \end{prooff}
                            
                            Similarly one can have another generalization when $K$ is connected.
                            \begin{thm}\label{eqcase}
                            Suppose that $\frak M$ is a neat collection
                            with only the lowest level possibly consisting of infinite components,
                            and that each component represents a nonzero class in the $\mathbb R$-homology of $X$. 
                            Let $K$ be a  connected compact  Lie transformation group of $X$.
                            Assume $\frak M$ is invariant under the action of $K$.
                            Then for any $K$-invariant Riemannian metric $g^K$, 
                            there exists a $K$-invariant metric $\hat g^K$ conformal to $g^K$
                            under which each component of  $\mathfrak M$ is homologically mass-minimizing.
                            \end{thm}
                            
            {\ }
                        
\subsection{More results}\label{mr}
              Since only one calibration is constructed for each dimension,
              results in \S\ref{sr}, e.g. Theorem \ref{many},
              lack the control on some region of the space of homology classes.
              To conquer this,
              we shall construct a metric that supports enough calibrations we need.

     When $X^n$ is oriented with betti number $b_k<\infty$ for $1\leq k<\frac{1}{2}n$,
     by Thom \cite{T} or Corollary II.30 in \cite{Thom}
     there exist embedded oriented connected compact  $k$-dimensional submanifolds
     $\mathcal L_k\triangleq\{M^k_1,\cdots,M^k_{b_k}\}$
     such that $span\{[M^k_i]\}_{i=1}^{b_k}=H_k(X;\mathbb R)$.
     By dimension reason one can arrange $\bigcup_{1\leq k <\frac{1}{2}n}\mathcal L_k$ to be a mutually disjoint collection.
      
       \begin{thm}\label{1}
       Let ${M^k_i}$ be given as above.
       Then in any conformal class of metrics
       there exists $\hat g$ under which
      every nonzero $\sum_{i=1}^{b_k} t_i[[M^k_i]]$ 
      where $1\leq k<\frac{1}{2}n$, $M^k_i\in\mathcal L_k$ and $t_i\in \mathbb R$ 
      is the unique mass-minimizing current in $\sum_{i=1}^{b_k} t_i[M^k_i]$.
      \end{thm}
       
       \begin{prooff}
       For the sake of simplicity, assume $\dim H_k(X;\mathbb R)$ $=$ $2$ for some $k< \frac{1}{2}n$ 
       and $\{[M_1],[M_2]\}$ is a basis where $M_1$ and $M_2$ are disjoint oriented connected compact submanifolds.
       Then there exist $k$-forms $\phi_1\text{ and }\phi_2$ on $X$ with $\int_{M_i}\phi_j=\delta_{i,j}.$
       Without loss of generality, assume $\phi_1\equiv 0$ on $U_{\epsilon}(M_2)$ and $\phi_2\equiv 0$ on $U_{\epsilon}(M_1)$ for some small $\epsilon$.
       Note that $\alpha$ can be chosen so that  
       $\|\hat \Phi_i\|^*_{\hat g}<\frac{1}{2}\text{ on } (U_\epsilon(M_i))^c$
       for the resulting forms $\hat \Phi_1$ and $\hat \Phi_2$ in \S\ref{forms} under the metric $\hat g$ in \S\ref{metrics}.
       A key observation is that
       $\pm\hat \Phi_1,\ \pm\hat \Phi_2\text{ and }
       \pm \hat \Phi_1\pm \hat \Phi_2$
       are all calibrations with respect to $\hat g$.
\begin{center}
\includegraphics[scale=0.28]{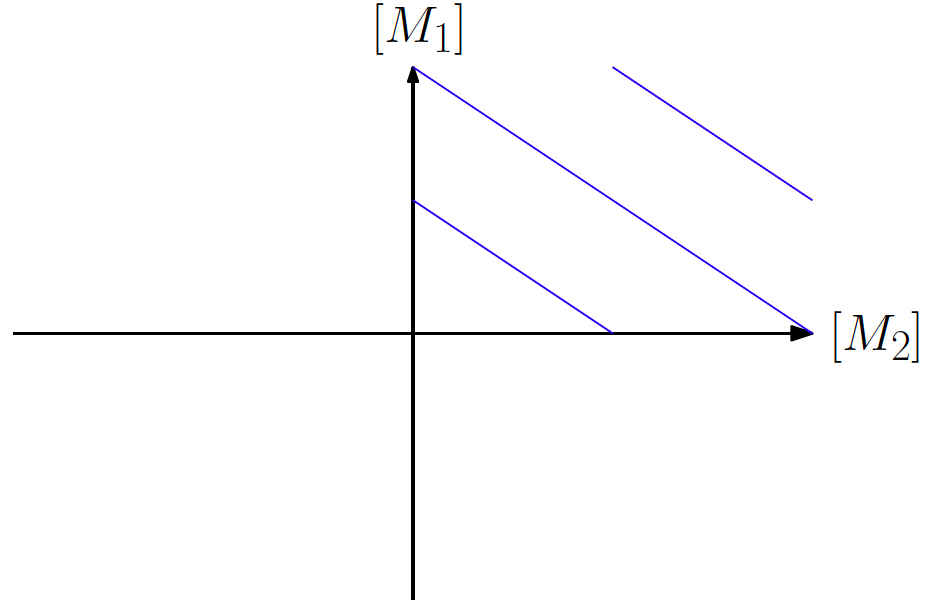}
\end{center}
        
        Then any nonzero linear combination of $[[M_1]]$ and $[[M_2]]$ can be calibrated in $(X,\hat g)$.
        For example, those representing classes in the closer of the first quadrant can be calibrated by $\hat \Phi_1+\hat\Phi_2$.
        The uniqueness follows as a result of 
$\mathrm{\bold{spt}}(\|\pm \hat \Phi_i\|^*_{\hat g}-1)= M_i$,
$\mathrm{\bold{spt}}(\|\pm \hat \Phi_1\pm \hat \Phi_2\|^*_{\hat g}-1)=\bigcup M_i$,
the simpleness of $\pm \hat \Phi_i$ along $M_i$
and $\pm \hat \Phi_1\pm \hat \Phi_2$ along $M_1\bigcup M_2$,
and the connectedness of $M_i$.

        When $\dim H_k(X;\mathbb R)=s$,
        $2^s$ such calibrations, each of which has comass norm bounded above 
        by $\frac{1}{s}$ away from some neighborhood of corresponding submanifold, 
        can be constructed for our purpose.
        More generally, for different dimension levels,
        the above argument combined with the elimination trick on forms in the proof of Theorem \ref{Xnoncthm} proves the theorem.
\end{prooff}

          When $\dim (X)\geq 6$, 
          one can choose $b_k$ smooth $k$-dimensional submanifolds $\mathcal L_k=\{M_i^k\}_{i=1}^{b_k}$ for $k=1,2,\cdots, n-3$,
          such that span$\{[M_i^k]\}_{i=1}^{b_k}=H_k(X,\mathbb R)$ and
such that
intersections $\mathcal I$ among $\bigcup^{n-3}_{k=1} \mathcal L_k$ are all transversal.
Note that $\mathcal I$ has a natural stratification structure 
$\cdots\prec \mathcal I_2\prec\mathcal I_1=\bigcup^{n-3}_{k=1} \mathcal L_k$,
where
$\mathcal I_t$ is the set of intersections among $t$ representatives.

\begin{thm}\label{2}
Let $X^n$ be an oriented manifold
with betti numbers $b_k<\infty$
for $1\leq k\leq n-3$
and $\mathcal L_k$ given above.
Then there exists a metric
such that 
every nonzero
$\sum_{i=1}^{b_k} t_i[[M^k_i]]$ where $1\leq k\leq n-3$, $M^k_i\in\mathcal L_k$ and $t_i\in \mathbb R$  
is the unique mass-minimizing current in
$\sum_{i=1}^{b_k} t_i[M^k_i]$.
\end{thm}
\begin{prooff}
       One can build a metric $g$ on $X$ such that, for any element $S$ of $\mathcal I_t$ ($t\geq 2$),
       there exists some $2\epsilon$-{cubic} neighborhood of $S$
       with fibers (induced by $g$ as in \S\ref{forms})
       split pointwise along $S$ as the Riemannian product of 
       fibers of $2\epsilon$-{cubic} neighborhoods of $S$ in $H_S$ for all $H_S \in \mathcal I_{t-1}$ and $S\subseteq H_S$.

        Let us focus on all (connected parts of) deepest intersections.
        For simplicity,
        suppose we have only one connected deepest intersection $\Delta$ and $\Delta\in\mathcal I_3$.
        Namely
        $\Delta$ is the intersection of three submanifolds $M_1$, $M_2$ and $M_3$.
        Assume $2\epsilon$ is universal for $S\in \bigcup_{t\geq 2}\mathcal I_t$ under $g$ in the preceding paragraph.
        Denote the volume form of $M_3$ by $\omega_3$,
        the distance function to $M_3$ by $\bold{d_3}$, and
        the projection to nearest point on $M_3$ by $\pi_3$.
        $\omega_1, \bold{d_1}, \pi_1$ and $\omega_2, \bold{d_2}, \pi_2$ are similarly defined.
        Since $\omega_i=d\psi_i$ in the $\epsilon$-neighborhood of 
        $(M_i\bigcap M_{i+1})\bigcup (M_i\bigcap M_{i+2})$  in $M_i$ (subscripts in the sense of mod $3$),
        define $\Psi_i=d(\rho_i\psi_i)$ in the union of $\epsilon$-cubic neighborhoods of  $M_i\bigcap M_{i+1}$ and $M_i\bigcap M_{i+2}$.
        Here we identify the pullback of $\omega_i$ (and $\psi_i$) via $\pi_i$ with itself,
        and $\rho_i$ is a smooth increasing function in $\bold{d_i}$
        with value zero when $\bold{d_i}\leq \frac{1}{2}\epsilon$ and value one for $\frac{2}{3}\epsilon\leq \bold{d_i}\leq\epsilon$.
        The slash-shadow region and the backslash-shadow region are
        intersections of regions $\Gamma_2:\frac{1}{2}\epsilon\leq\bold{d_2}\leq\frac{2}{3}\epsilon$
        and $\Gamma_3:\frac{1}{2}\epsilon\leq\bold{d_3}\leq\frac{2}{3}\epsilon$ 
        with $M_1$ respectively.

\begin{center}
\includegraphics[scale=0.38]{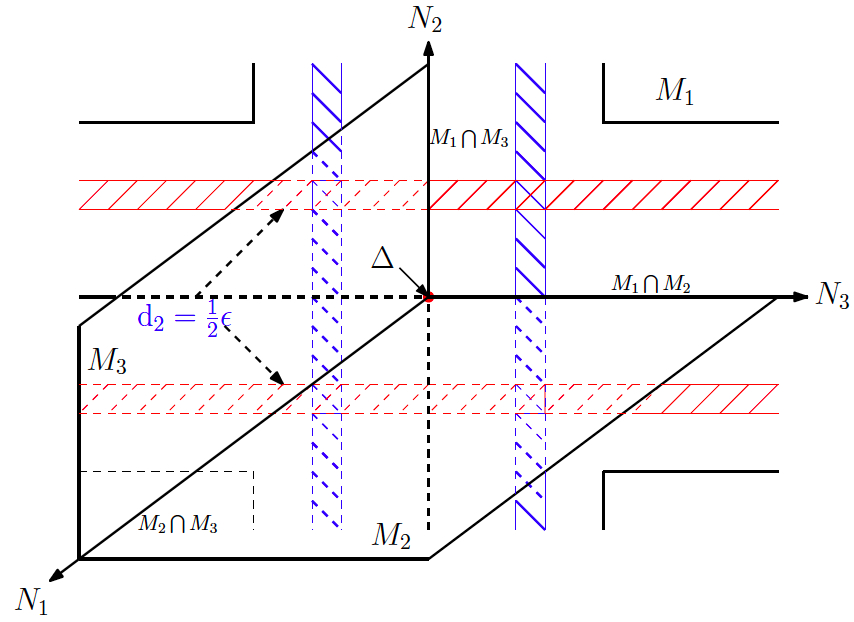}
\end{center}

         \textbf{Case One:} Every $M_i$ has the same dimension $k$.
         There are three (bunches of) directions $N_i$ on the $\epsilon$-cubic neighborhood $U_\epsilon(\Delta)$ 
        of $\Delta$.
        (Note that  $N_i$ has meanings for $\bold{d_i}\leq 2\epsilon$ only.)
        Denote the split part of $g$ along $N_i$ by $g_i$.
        
        \textbf{Claim:} 
                $\sum\omega_i$ is a calibration in $U_\epsilon(\Delta)$.
         
         Since a form and its Hodge dual have the same comass,
         the claim is an immediate consequence by applying the following lemma to $*\sum\omega_i=\sum*\omega_i$.
         
         \begin{lem}\label{sigmaomegacal}
         Let $e_1,\cdots, e_{n+2}$ be an orthonormal basis for $\mathbb R^{n+2}$,
         and for each multi-index $I=\{i_1,\cdots, i_p\}$ where $i_1<\cdots<i_p$,
         let $e_I^*$ denote the corresponding ``axis" $p$-form $e_{i_1}^*\wedge\cdots\wedge e_{i_p}^*$.
         Assume $\phi= e^*_J\wedge e^*_{n+1}\wedge e^*_{n+2}$ where $J=\{j_1,\cdots, j_{p-2}\}\subset \{1,\cdots, n\}$
         and $\psi=\sum_I e^*_I$ with $i_p\leq n$.
         Then $$\|\phi+\psi\|^*=\max \{1,\|\psi\|^*\}.$$
         \end{lem}       
         However $\sum\omega_i$ is not well defined on the union $\Xi$ of $\epsilon$-cubic neighborhoods of  $M_1\bigcap M_2$, $M_2\bigcap M_3$ and $M_3\bigcap M_1$.
         Instead we consider 
                \begin{equation*}
                \phi_k=\sum \omega_i-\sum \Psi_i=\sum\big[(1-\rho_i)\omega_i-d\rho_i\wedge\psi_i\big] \text{ on } \Xi.
                \end{equation*}
         Then 
                 \begin{equation}\label{phi}
                  \phi_k=\omega_i \text{ when } \bold{d_{i}}\leq\frac{1}{2}\epsilon, \text{ } \bold{d_{i+1}}\geq\frac{2}{3}\epsilon
                   \text{ and } \bold{d_{i+2}}\geq\frac{2}{3}\epsilon.
                    \end{equation}

        Note that, under the condition $n-k\geq 3$, for example in $M_1$,
        the subspace spanned by the dual $k$-vector of each simple form $(1-\rho_i)\omega_i-d\rho_i\wedge\psi_i$ for $i\neq1$
        contains at least 2 directions of $N_1$.
        So, by the useful result of Harvey and Lawson below,
        if one multiples $g_1$ by a sufficiently large constant $C>1$,
        then $\phi_k$ has comass one (same as that of $\omega_1$) in $\Xi\bigcap M_1$.
        
                    \begin{lem}[Corollary 2.11. in \cite{HL1}]\label{sumcomass}
                    With notation as in Lemma \ref{sigmaomegacal},
                   $$ \|e^*_1\wedge \cdots \wedge e^*_p+\sum_I b_I e^*_I\|^*\leq \max\{1,\sum_I |b_I|\}$$
                    provided that $b_I=0$ whenever $i_{p-1}\leq p$.
                    \end{lem}

         We are now about to modify $g$ so that $\phi_k$ becomes a calibration in some neighborhood of $\bigcup M_i$.      
         Let $C$ work for every $\Xi\bigcap M_i$.       
        Choose a smooth function $f$ of $\bold{d}$ for $\bold{d}\leq\epsilon$ as in the picture and set $f_i=f(\bold{d_i})$.
\begin{center}
\includegraphics[scale=0.4]{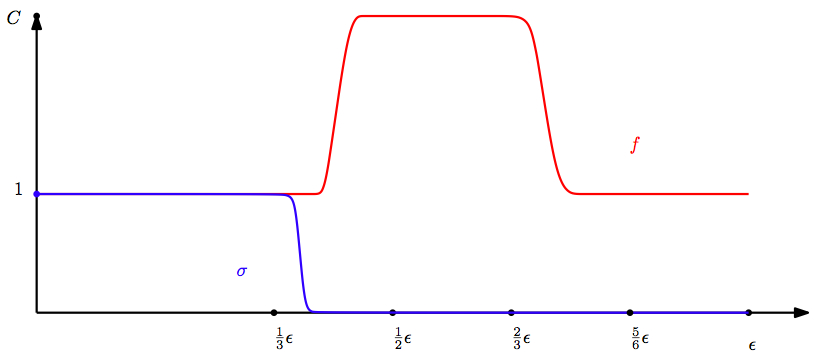}
\end{center}        
                  Along $M_1\bigcap M_2$,
           set
           $$g_1\rightarrow f_3g_1,\quad g_2\rightarrow f_3g_2,\quad 
           \text{ and } \quad g_3\rightarrow f_3^{-1}g_3\quad (\star)$$
           and similarly for $M_2\bigcap M_3$ and $M_3\bigcap M_1$.
           Then in these three sets $\phi_k$ becomes a calibration.
           
           We want to extend the metric along each $M_i$.
           A good try based on $(\star)$ to the $\epsilon$-neighborhood of $M_1\bigcap M_2$ in $M_1$ is this.
             $$g_1\rightarrow f_3g_1,\quad g_2\rightarrow f_3^{\sigma_2}g_2,\quad 
           \text{ and } \quad g_3\rightarrow f_3^{-{\sigma_2}}g_3\quad (\ast)$$
           where $\sigma$ is a cutoff function with $f=1$ on $\bold{spt}(\sigma)$ and $\sigma_i=\sigma(\bold{d_i})$.
           A subtle point here is that the volume form of $M_1$ is unchanged.
           The same extension from $M_1\bigcap M_3$ to $M_1$ gives the following.
            $$g_1\rightarrow f_2g_1,\quad g_2\rightarrow f_2^{-\sigma_3}g_2,\quad 
           \text{ and } \quad g_3\rightarrow f_2^{{\sigma_3}}g_3.\quad $$
           Since these two extension do not agree in $U_\epsilon(\Delta)\bigcap M_1$, we combine them by
           $$g_1\rightarrow f_2f_3g_1,\quad g_2\rightarrow f_2^{-\sigma_3}f_3^{\sigma_2}g_2,\quad 
           \text{ and } \quad g_3\rightarrow f_2^{{\sigma_3}}f_3^{-{\sigma_2}}g_3.\quad $$

           By the choice of $C$, $f\geq 1$ and  $f=1$ on $\bold{spt}(\sigma)$,
           the above combination makes the comass of $\phi_k$ one under the modified metric following Lemma \ref{sumcomass}. 
           Together with the same procedure for $M_2$ and $M_3$, we get a metric $\tilde g$ in $\bigcup\Xi\bigcap M_i$
           that makes $\phi_k$ a calibration.
           Note that in the intersection of $\frac{1}{3}\epsilon$ neighborhoods of $M_{i+1}$ and $M_{i+2}$ in $M_i$,
           or in the complement of the intersection of $\frac{5}{6}\epsilon$ neighborhoods of $M_{i+1}$ and $M_{i+2}$ in $U_\epsilon(\Delta)\bigcap M_i$,
           $\tilde g=g$.
           Hence by \eqref{phi} $\tilde g$ 
            produces a metric $\check g$
           on the union $\Upsilon$ of $\frac{1}{3}\epsilon$ neighborhoods of $M_i$ (containing each $M_i$)
           making $\phi_k$ a calibration.
           Furthermore every nonzero $\sum_{i, n_i=\pm 1,0}n_i(\omega_i-\Psi_i)$ becomes a calibration in $(\Upsilon, \check g)$.
\\{\ }
 
        \textbf{Case Two:} 
        $M_2$ and $M_3$ are of dimension $k$, but $M_1$ has a different dimension $m$.
        (Similar for the case with mutually different dimensions.)
        Consider potential calibrations 
        $ \pm(\omega_2-\Psi_2), \pm(\omega_3-\Psi_3),\pm(\omega_2-\Psi_2)\pm(\omega_3-\Psi_3), \text{ and }\pm(\omega_1-\Psi_1)$ on $\Xi$.
        In a similarly way with different weights in $(\star)$ and $(\ast)$
        one can get calibration pairs on some neighborhood of $\bigcup M_i$.
\\{\ }

     The idea works for general cases with modified $(\star)$ and $(\ast)$.
     Following the above steps around all connected parts of deepest intersection,
     one can extend the preferred local calibration pairs to global ones sharing a common metric.
     Multiply the metric by a smooth function which is one in $\bigcup_{k=1}^{n-3}\bigcup_{M\in\mathcal L_k}M$ and strictly greater than one elsewhere.
     Name it $\hat g$.
     Then every nonzero $\sum_{i=1}^{b_k} t_i[[M_i^{k}]]$ for  $1\leq k\leq n-3$ and $t_i\in \mathbb R$ can be calibrated in $(X,\hat g)$.
     The uniqueness of such a mass-minimizing current in its current homology class follows similarly as in the proof of Theorem \ref{1}.
     Here note that for any point $p\in M^k_i-\mathcal I_2$ (a.e. on $M^k_i$)
     the oriented unit $k$-vector of $\wedge^k{T_pM_i^k}$ is the unique unit $k$-vectors in $\wedge^kT_pX$
     that has pairing value one with the corresponding calibrations of $M_i^k$.
\end{prooff}
       \begin{rem}
       To illustrate how $(\star)$ and $(\ast)$ change,
       suppose in Case One we have $M_{123}$, $M_{124}$, $M_{134}$, $M_{234}$ and perpendicular directions $N_4$, $N_3$, $N_2$, $N_1$ respectively given in the figure below.
\begin{center}
\includegraphics[scale=0.5]{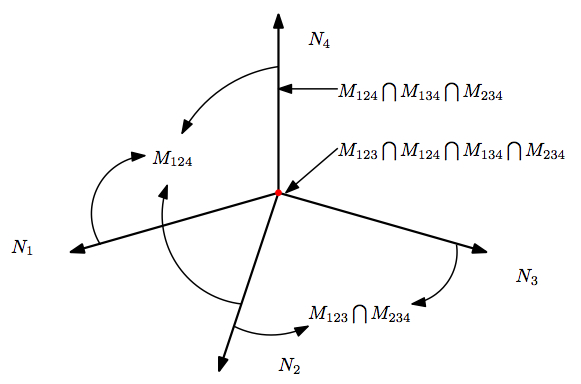}
\end{center}         
       Then, along $M_{124}\bigcap M_{134} \bigcap M_{234}$, $(\star)$ transforms to
             $$g_1\rightarrow f_4g_1,\quad g_2\rightarrow f_4g_2,\quad g_3\rightarrow f_4g_3 \quad
           \text{ and } \quad g_4\rightarrow f_4^{-2}g_4\quad (\star')$$
           and a good try of metric extension to $M_{124}$ is
             $$g_1\rightarrow f_4^{\sigma_1\sigma_2}g_1,\quad g_2\rightarrow f_4^{\sigma_1\sigma_2}g_2,\quad g_3\rightarrow f_4g_3 \quad
           \text{ and } \quad g_4\rightarrow f_4^{-2\sigma_1\sigma_2}g_4\quad (\ast').$$
       The corresponding $\tilde g$ in $U_\epsilon(\Delta)\bigcap M_{124}$ is given by
\begin{equation*}     
             \begin{split}
                                        &g_1\rightarrow f_1^{-2\sigma_2\sigma_4}f_2^{\sigma_1\sigma_4}f_4^{\sigma_1\sigma_2}g_1,\\
                                                                &g_2\rightarrow f_1^{\sigma_2\sigma_4}f_2^{-2\sigma_1\sigma_4}f_4^{\sigma_1\sigma_2}g_2,\\
                                                                &g_3\rightarrow f_1f_2f_4g_3, \text{ and }\\
                                                                &g_4\rightarrow f_1^{\sigma_2\sigma_4}f_2^{\sigma_1\sigma_4}f_4^{-2\sigma_1\sigma_2}g_4.
                                     \end{split}
\end{equation*}
                   \end{rem}
       
        
   \begin{rem}
     Generally codimension at least 3 is necessary for applying Lemma \ref{sumcomass}.
     When $n=4$ or $5$, Theorem \ref{2} can be improved to include the level of codimension $2$ by Theorem \ref{conecal}.
     \end{rem}
     
     \begin{pf6}
     Assume the comass of $\phi+\psi$ is achieved by
     pairing with a unit $p$-vectors $\xi$.
     Then we will make use of the following ``canonical form of a simple vector with respect to a subspace".
     \begin{lem}[Lemma 7.5 in \cite{HL2}]\label{cfsv}
     Suppose $V\subset \mathbb R^n$ is a linear subspace and $\xi$ is a unit simple $p$-vector.
     Then there exists set of orthonormal vectors $f_1,\cdots,f_r$ in $V$, a set of orthonormal vectors $g_1,\cdots, g_s$ in $V^\perp$,
     and angles $0<\theta_j<\frac{\pi}{2}$ for $j=1,\cdots, k$ (where $k\leq r,s\leq p$ and $r+s-k=p$)
     such that
     $$\xi=(\cos \theta_1f_1+\sin\theta_1 g_1)\wedge\cdots\wedge (\cos \theta_kf_k+\sin\theta_k g_k)
     \wedge f_{k+1}\wedge\cdots\wedge f_r\wedge g_{k+1}\wedge\cdots\wedge g_s.$$
     \end{lem}
     {\ }\\
     Let $V=span\{e^{n+1},e^{n+2}\}$.
     These $\lambda_j=\cos^2\theta_j$ are eigenvaules of a symmetric bilinear form $B$
     where $\pi:\mathbb R^m \rightarrow V$
     and $B(u,v)=<\pi(u),\pi(v)>$ is defined on span $\xi$. 

            Assume $r=k=2 \text{ and }s=p$ (otherwise either $<\phi,\xi>$ or $<\psi, \xi>$ gives zero and a proof or contradiction simply follows).
            We have
            $$\xi=(\cos \theta_1e_1+\sin\theta_1 g_1)\wedge(\cos \theta_2e_2+\sin\theta_2 g_2)
     \wedge g_{3}\wedge\cdots\wedge g_p.$$
          Evaluating $\phi+\psi$ on $\xi$ leads to
\begin{equation*}     
             \begin{split}
     |<\phi+\psi,\xi>|&=|\cos\theta_1\cos\theta_2\cdot\phi(e_1,e_2,g_3\cdots,g_p)+\sin\theta_1\sin\theta_2\cdot\psi(g_1,g_2,g_3\cdots,g_p)|\\
                             &\leq \cos\theta_1\cos\theta_2\cdot|\phi(e_1,e_2,g_3\cdots,g_p)| +\sin\theta_1\sin\theta_2\cdot|\psi(g_1,g_2,g_3\cdots,g_p)|\\
                             &\leq \cos(\theta_1-\theta_2)\cdot \max\{1,\|\psi\|^*\}\\
                             &\leq \max\{1,\|\psi\|^*\}
   \end{split}
\end{equation*}
     \end{pf6}
\textbf{Question:}
Usually one cannot have such existence result when $k$ can be $n-1$.
Therefore it may be interesting to ask whether the same conclusion holds for $1\leq k\leq n-2$ in general.  
\\{\ }

                 \section{Singular case}\label{SS}
                 In this section the case of submanifolds with singularities will be discussed.
                 Unlike the smooth case, one cannot have local calibration pairs so easily as in \S\ref{poc}.
                 Our concern here is to extend an existing local calibration pair around the singular set
                 to a calibration pair on some neighborhood of the singular submanifold under consideration.
                 Then a further extension from the neighborhood to global is roughly the same as in the smooth case.
                 
                 We first recall two useful lemmas by Harvey and Lawson, then obtain our extension theorem, and next apply
                 it for several interesting examples in the realm of calibrated geometry.
                         
         
\subsection{Two lemmas}\label{tl} %
            The first lemma tells us how to canonically decompose a $p$-form with respect to certain $p$-plane. 
            \begin{lem}[Lemma 2.12 in \cite{HL1}] \label{hl1}
             Let $\xi\in \Lambda^p \mathbb{R}^n$ be a simple p-vector with
             $V = span\ \xi$. Suppose $\phi\in\Lambda^p\mathbb{R}^n$ satisfies  $\phi(\xi)= 1$. 
              Then there exists a unique oriented complementary subspace $W$ to $V$  with the following property.
                      For any basis $v_1, \cdots, v_n$ of $\mathbb{R}^n$ such that $\xi=v_1\wedge... \wedge v_p$
                      and $v-{p+1}, \cdots, v_n$ is basis for $W$, 
                      one has that
                                 \begin{equation}
                                 \phi=v_1^*\wedge \cdots \wedge v_p^*+ \sum a_Iv_I^*,
                                  \end{equation}
                                   where $a_I=0$ whenever $i_{p-1}\leq p$. 
                                   Here $I=\{i_1,\cdots, i_p\}$ with $i_1<\cdots<i_p$.
              \end{lem}
             The second lemma says how to create metrics based on the above decomposition with control on the comass of the form. 
             A combination of Lemma \ref{hl1} and Lemma \ref{sumcomass} gives its proof.
             \begin{lem}[Lemma 2.14 in \cite{HL1}]\label{hl2}
                    Let $\phi,\  V=span \ \xi$, and $W$ be given above.
                    Consider an inner product $<\cdot, \cdot>$ on $\mathbb{R}^n$ such that $V\perp W$ and $\|\xi\|=1$. 
                    Choose any constant $C^2 > \bigl( \begin{smallmatrix} n\\ p\end{smallmatrix} \bigr) \|\phi\|^*$
                    and define a new inner product on $\mathbb{R}^n=V\oplus W$
                    by setting $<\cdot,\cdot>'=<\cdot,\cdot>_V+C^2<\cdot,\cdot>_W$.
                    Then under this new metric we have
                             $$\|\phi\|^*=1\ and\ \phi(\xi)=\|\xi\|=1.$$
             \end{lem}
             \begin{rem}\label{imp}
              If 
 $\phi(\xi)=\vartheta$ (positive) not necessarily one,
 one can apply {\it Lemma \ref{hl1}} to $\vartheta^{-1}\phi$ for
 $\|\phi\|^*=\vartheta,\ \|\xi\|=1\ and\ \phi(\xi)=\vartheta$
 by choosing $C^2>\vartheta^{-1} \bigl( \begin{smallmatrix} n\\ p\end{smallmatrix} \bigr)
\|\phi\|^*$.
              \end{rem}
              
              They will be used in proving the extension result in the next subsection.

{\ }


\subsection{An extension result}\label{CWS}%
        \begin{defn}       
        By a \textbf{singular\ submanifold} $(S,\mathscr S)$ of dimension $m$ 
        with singular set $\mathscr S$, we mean a pair of closed subsets $\mathscr S\subset S$ of X,
        where $S-\mathscr S$ is an $m$-dimensional submanifold and the Hausdorff $m$-measure $\mathcal H^m(\mathscr S) = 0$.
        We say $(S,\mathscr S)$ is \textbf{orientable}, if $S-\mathscr S$ is so.
        \end{defn}

          \begin{rem}\label{coneccc}
          Assume $S$ is an oriented compact submanifold with only one singular point $p$ and $C_p$ is a tangent cone of $S$ at $p$.
          Then the current $[[S]]=\int_S\cdot\ $ is calibrated by a smooth $\phi$ if and only if $S-p$ is calibrated by $\phi$.
          Moreover, either of them implies that $\phi_p$ calibrates $C_p$ in $(T_pX,g_p)$.
          \end{rem}
          
          From now on, $(S,o)$ will be assumed an oriented connected compact singular submanifold with one singular point $o$.
                 \begin{figure}[ht]
                 \begin{center}
                 \includegraphics[scale=0.26]{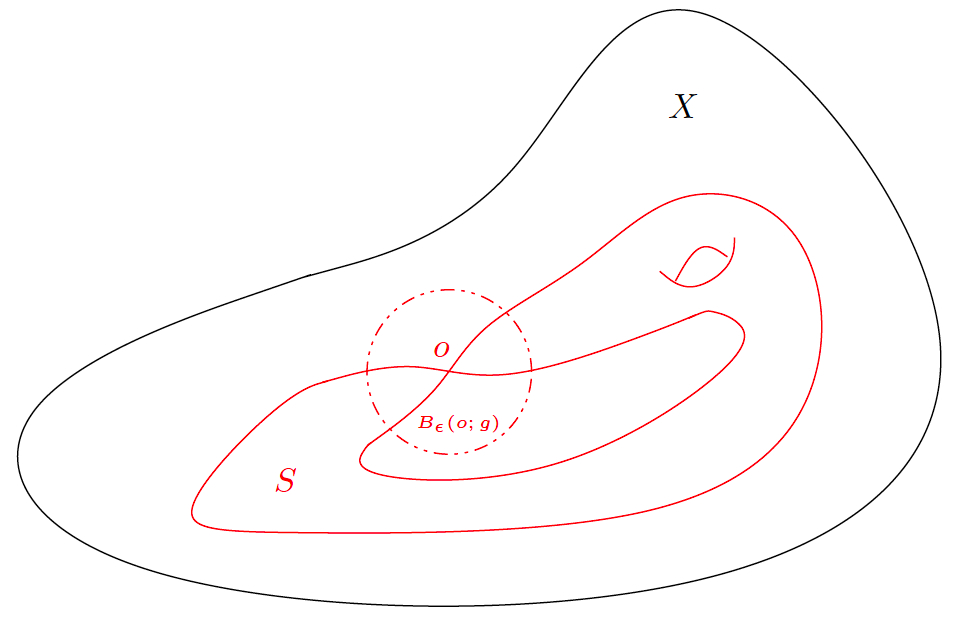}\ \ \ \ \ \ \ 
                 \end{center}
                 \end{figure}
                            
          \begin{thm}\label{conecal}
          Given $(S,o)$ in $(X,g)$
          with $[S]\neq[0]\in H_m(X; \mathbb{R})$.
           If  $B_{\epsilon}(o;g)\cap S$
           can be calibrated by a smooth calibration in 
            some $\epsilon$-ball $(B_\epsilon(o;g),g)$ centered at $o$,
           then
            there exists a metric $\hat g$ coinciding with $g$ on $B_{{\frac{\epsilon}{2}}}(o;g)$ such that
            $S$ can be calibrated by a smooth calibration in $(X,\hat g)$.
            \end{thm}

       \begin{rem}
       In the theorem, 
       $\frac{\epsilon}{2}$ can be replaced by $\varkappa\epsilon$ for any $0<\varkappa<1$.
       \end{rem}

               \begin{prooff}
               Assume $\epsilon$ is small enough so that 
               the local calibration $\phi$ on $B_\epsilon(o;g)$ 
               can be written as $d\psi$ for some smooth $(m-1)$-form $\psi$.
                Suppose the compact region $\Gamma_1\bigcup \Omega\bigcup\Gamma_2$
                (the diffeomorphic image of an $h$-disk normal bundle, for small $h$, over
                a closed set $(\Gamma_1\bigcup \Omega\bigcup\Gamma_2)\cap S$ 
                by the exponential map restricted to normal directions, see picture below)
                is contained in $B_\epsilon(o;g)-B_{{\frac{2\epsilon}{3}}}(o;g)$.
 \begin{figure}[h]
\begin{center}
\includegraphics[scale=0.3]{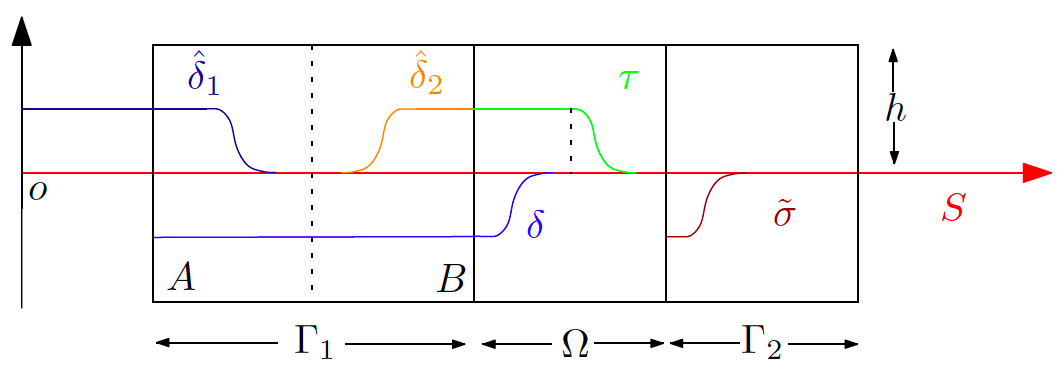}
\end{center}
\end{figure}
           Denote the projection by $\pi$ and call the directions perpendicular to fibers {horizontal}.
           Then $\pi^*\omega= d(\pi^*(\psi|_S))$ in $\Gamma_1\bigcup \Omega\bigcup\Gamma_2$
            where $\omega$ is the volume form of $S\bigcap(\Gamma_1\bigcup \Omega\bigcup\Gamma_2)$. 
            Set
 \[
 \Phi\triangleq d(\tau\psi+(1-\tau)\pi^*(\psi|_S))=d\tau\wedge(\psi-\pi^*(\psi|_S))+\tau\phi+(1-\tau)\pi^*\omega
 \]
          where $\tau$ is a cutoff function shown in the picture
          with value one near $\Gamma_1$ and zero near $\Gamma_2$.
          (The picture is just an illustration, since ``height" $h$ is usually smaller than one.)
 
           Denote by $\overrightarrow{H_{y,g}}$ the unique oriented unit horizontal $m$-vector at $y$ with respect to $g$
           for $y \in \Gamma_1\bigcup \Omega\bigcup\Gamma_2$.
           Observe that $\Phi(\overrightarrow{H_{x,g}})=1$ for $x\in S\bigcap(\Gamma_1\bigcup \Omega\bigcup\Gamma_2)$.
             By shrinking $h$,
             the smooth function $\Phi(\overrightarrow{H_{y,g}})>\frac{1}{2}\text{ on }\Gamma_1\bigcup \Omega\bigcup\Gamma_2$.
 Set 
 \[
 \bar g=f\cdot g
\text{ where } 
f=\delta+(1-\delta)(\Phi(\overrightarrow{H_{y,g}}))^{\frac{2}{m}}
\]
          on $\Gamma_1\bigcup \Omega\bigcup\Gamma_2$.
          Note $f=1$ in $S\bigcap(\Gamma_1\bigcup \Omega\bigcup\Gamma_2)$ and $\Phi=\phi$ on $\bold{spt}(\delta)$.
          Since $(\phi,g)$ is a local calibration pair,
          $f\geq (\Phi(\overrightarrow{H_{y,g}}))^{\frac{2}{m}}$ in $\Gamma_1\bigcup\Omega\bigcup \Gamma_2$ and $f\equiv1$ in $\Gamma_1$.
          Then $\Phi$ and $\bar g$ naturally extend on $\Upsilon$, the region embraced by the ``curve" in the picture below
          (an ``$h$-disk bundle'' containing $\Gamma_1\bigcup\Omega\bigcup\Gamma_2$).
           \begin{figure}[h]
\begin{center}
\includegraphics[scale=0.26]{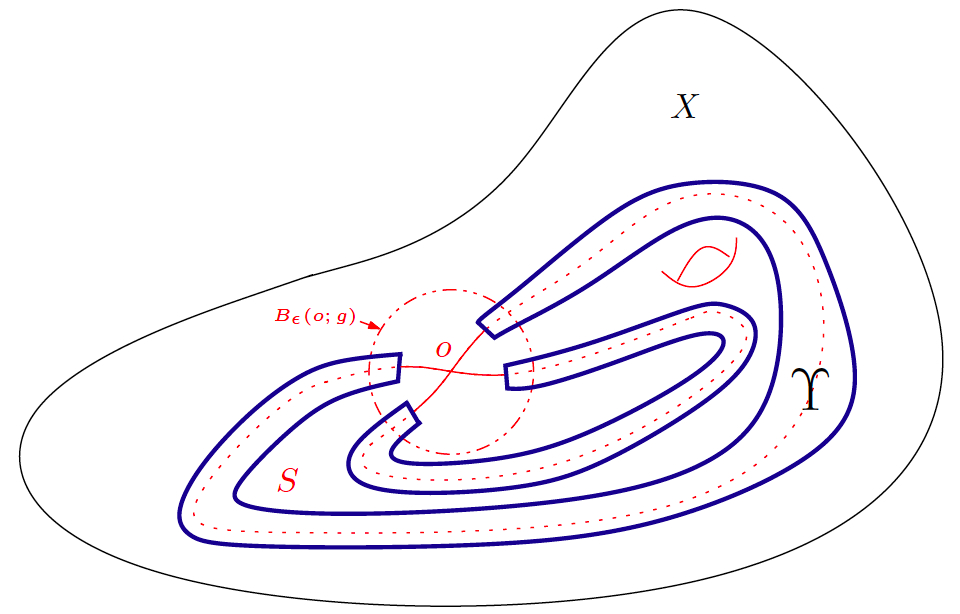}
\end{center}
\end{figure}
Note that

(a). $\Phi$ calibrates $S\bigcap(\Upsilon-\Omega)$ in $(\Upsilon-\Omega,\bar g)$,

(b). $\bar g=g$ in $\Gamma_1$, and

(c). $\frac{1}{2}<\Phi(\overrightarrow{H_{y,{\bar g}}})\leq 1$ on $\Upsilon$ with equality on $\Upsilon-\Gamma_1-\Omega$,
where $\overrightarrow{H_{y,{\bar g}}}$ is the unique oriented unit horizontal $m$-vector at $y$ for $\bar g$.

{\ }
 
      Now we wish to glue $\bar g$ and $g$ together to make $\Phi$ a calibration.
      By applying {Lemma \ref{hl1}} to $\Phi$, $\overrightarrow{H_{y,{\bar g}}}$ and $\bar g$ on $\Upsilon$,
      there is a smoothly varying $(n-m)$-dimensional plane field $\mathscr W$ transverse to horizontal directions in $\Upsilon$.
      Following {Lemma \ref{hl2}}, {Remark \ref{imp}} and Property (c), for any metric $g_{\mathscr W}$ along $\mathscr W$,
      there exists a sufficiently large constant $\bar \alpha$ (due to the compactness of $\Upsilon$) such that,
      under $\tilde g= \bar g^h\oplus\bar \alpha g_{\mathscr W}$ on $\Upsilon$, 
         where $\bar g^h$ is the horizontal part of $\bar g$,
\[
 \|\Phi\|_{y,\tilde g}^*=
\Phi(\overrightarrow{H_{y,{\bar g}}})\leq 1\text{ for } y\in \Upsilon,\text{ and } \Phi(\overrightarrow{H_{x,{\bar g}}})=1\text{ for } x\in S\bigcap \Upsilon.
\]

          \begin{figure}[h]
          \begin{center}
          \includegraphics[scale=0.26]{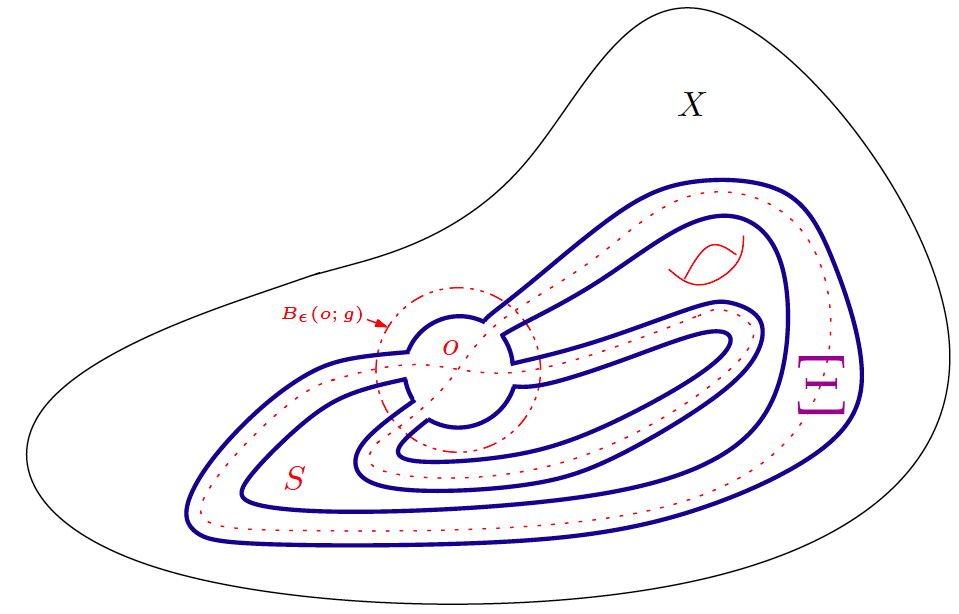}
          \end{center}
          \end{figure}
         
         Based on Property (b) we construct a smooth metric $\check g$ on $\Xi$ as follows.
         $$
         \check g=
              \begin{cases}\label{aaa}
               g & \text{near}\ o\\
               g+(1-\hat \delta_1)((0\cdot {\bar g^h})\oplus \bar \alpha g_{\mathscr W}) &
               {\text{on} \ A-\Gamma_0}\\
        {       (0\cdot g^h)\oplus g^\nu+ (0\cdot \bar g^h)\oplus \bar \alpha g_{\mathscr W}
               +\hat\rho (g^h\oplus (0\cdot g^\nu))+(1-\hat \rho)(\bar g^h\oplus (0\cdot g_{\mathscr W}))}
               & 
                 {\text{on } \Gamma_0}\\
               (1-\hat \delta_2)((0\cdot g^h)\oplus g^\nu) +\tilde g &
               {\text{on} \ B-\Gamma_0}\\
               \tilde g & \text{on}\ \Omega\\
               \tilde \sigma
               \tilde g+(1-\tilde \sigma)
               \bar g & \text{on}\ \Gamma_2\\
               \bar g & \text{far away from}\ o\\
        \end{cases}
         $$
           Here
          {$\hat \rho$ is a deceasing function from value $1$ to $0$ with $supp(d\hat\rho)\Subset$$\Gamma_0\triangleq\{\hat \delta_1=0\}\bigcup \{\hat \delta_2=0\}$,}
         $g^h,g^\nu$ are the horizontal and fiberwise parts of $g$ respectively,
          $\oplus$ means the orthogonal splitting of a (pseudo-)metric and $+$ is the usual addition between two (pseudo-)metrics.  
          Note that, 
      on $\Gamma_2$,  ${\mathscr W}$ is exactly the distribution of fiber directions
      and $\Phi=\pi^*(\omega)$ is a simple horizontal $m$-form;
      {
      and that $\mathscr W$ is also vertical in $\Gamma_0\bigcap S$.
          Shrink $h$ for $\|\phi\|^*_{\check g}$ being smooth in $\Gamma_0$
          and replace $\check g$ by $[\hat \delta (\|\phi\|^*_{\check g})^{\frac{2}{m}}+(1-\hat\delta)]\check g$ on $\Gamma_0$
          where $\hat \delta$ is a bump function of value $1$ on $supp(d\hat\rho)$ and $0$ near the ends of $\Gamma_0$.
          Now $\Phi$ becomes a calibration in $\Xi$.}
      Since $S$ is a strong deform retract of $\Xi$ and $[S]\neq 0$,
      it can extend to a global calibration pair of $S$ by \S\ref{forms} and \S\ref{metrics}.
\end{prooff}

      By observing that the comass function of a smooth form of co-degree one is always smooth, we have the following refinement.
            \begin{cor}\label{codim1}
            Suppose $(S,o)$ is of codimension one in $(X,g)$ representing a nonzero real homology class.
            If $B_{\epsilon}(o;g)\cap S$ for some $\epsilon>0$ can be calibrated by a coflat calibration singular only at $o$ in $(B_\epsilon(o;g),g)$, 
            then there exists a metric $\hat g$ conformal to $g$ with $\hat g=g$ on $B_{{\frac{\epsilon}{2}}}(o;g)$
            such that $S$ can be calibrated by a coflat calibration singular only at $o$ in $(X,\hat g)$.
            \end{cor}

            In fact it does not have to require that $S$ is a strong deformation retract of some open neighborhood of $S$ for the last step in the proof. 
                  \begin{cor}
                  Suppose $(S,\mathscr S)$ is of dimension $m$ in $(X,g)$.
                  Assume $V\bigcap S$
                  for some open neighborhood $V$ of $\mathscr S$
                  can be calibrated in $(V,g|_V)$ by some coflat calibration $\phi$ with singular set $S_\phi\subset \mathscr S$.
                  Assume further $[S]\neq[0]$ in $H_m(X;\mathbb R)$.
                  If 
                        \begin{equation}\label{i}
                        i^*: H^m(X-\mathscr S;\mathbb R)\rightarrow H^m(U_\epsilon-\mathscr S;\mathbb R),
                        \end{equation}
                  is surjective for every $\epsilon$-neighborhood $U_\epsilon$ of $S$ under $g$ where $0<\epsilon\leq\epsilon_0$ for some small $\epsilon_0$,   
                   then there exists a metric $\hat g$
                    such that
                     $S$ can be calibrated in $(X,\hat g)$ by a coflat calibration with singular set $S_\phi$.
                  \end{cor}
               
                \begin{rem}
                By {\it Almgren}'s big regularity theorem, 
                being calibrated of $S$ around $\mathscr S$ implies that
                $\mathscr S$ has codimension at least $2$ in $S$.
                By $\bold{spt}(d[[S]])\subseteq\mathscr S$, $d[[S]]=0$ automatically and therefore $[S]$ makes sense.
                \end{rem}

                \begin{rem}
                When $\mathscr S$ is a smooth submanifold, 
                $S$ is a strong deformation retract of $U_{\epsilon}$ for small $\epsilon$.
\end{rem}

               \begin{prooff}
                     Choose $\delta>0$ such that $5\delta$-neighborhood of $\mathscr S$ is contained in $V$
                     and $5\delta\leq \epsilon_0$.          
                     Let $f$ be a $\delta$-approximation Morse function of the distant function to $\mathscr S$ in $S-\mathscr S$ in $C^0$ norm. 
                     Then, on some ``annulus-like" region 
                     $ \{    \delta<a_1\leq f\leq a_2<4\delta    \}$
                     which contains no critical points,
                     the volume form of $S$ is exact.
                     By a suitable choice of gluing function $\tau(f)$,
                     an extension of $(\phi, g|_V)$ to a coflat calibration pair about the entire $(S,\mathscr S)$ can be made.
                     Furthermore, by the subjectivity of \eqref{i}, a global calibration pair can be produced without new singularities.
               \end{prooff}
{\ }

\subsection{Further applications}\label{exs}
           Under some circumstances calibrations cannot avoid having singularities.
           In \cite{Z4} we showed that
           every homogeneous area-minimizing hypercones
           can have calibrations singular only at the origin.
\\

\textbf{Example 1:} 
          When the local model around $o$ in Theorem \ref{conecal} is a {Simons} cone over $S^{r-1}\times S^{r-1}$ for $r\geq 4$,
          one has a smooth calibration $ \phi$ (which actually can be $SO(r)\times SO(r)$ invariant) on $\mathbb R^{2r}-\{0\}$.
          Follow the proof of Theorem \ref{conecal} to get $\Phi$ on $\Xi-o$ and $\check g$ on $\Xi$. 
          The {\em Mayer-Vietoris} sequence for $\Xi-o$ and an open ball $B$ with $o\in B\subset \Xi$
          gives the following the exact sequence
 $$\cdots \rightarrow H^{2r-2}(S^{2r-1}(\upsilon))\rightarrow H^{2r-1}(\Xi)\rightarrow H^{2r-1}(\Xi-o)\rightarrow H^{2r-1}(S^{2r-1}(\upsilon))\rightarrow H^{2r}(\Xi)\rightarrow \cdots$$
          where $S^{2r-1}(\upsilon)$ is a small $\upsilon$-sphere centered at $o$.
          Since
                  \[
                  \|\int_{S^{2r-1}(\upsilon)}\Phi\|=
                  \|\int_{S^{2r-1}(\upsilon)}\phi\|=
                  \lim_{\upsilon\rightarrow 0}\|\int_{S^{2r-1}(\upsilon)}\phi\|     \leq 
                  \lim_{\upsilon\rightarrow 0}\text{ vol}(S^{2r-1}(\upsilon))=0,
                  \]
            $S$ is a strong deformation retraction of $\Xi$, and $[S]\neq 0$,
             there is a smooth form $\check \phi$ on $X$ such that
                \[
                 \check \phi|_{\Xi-o}-\Phi=d\check \psi
                \]
                 for some smooth $(2r-2)$-form $\check \psi$ on $\Xi-o$.
           Now, away from $S$, glue $\check \phi$ and $\Phi$ together to a smooth form $\hat \Phi$ on $X-o$,
           and meanwhile extend $\check g$ to $\hat g$ making $\hat\Phi$ a calibration on $X-o$.

         By Remark \ref{coflatFTCG}, $[[S]]$ is homologically mass-minimizing.
         However, it is impossible to calibrate $S$ using a smooth calibration $\bar \Phi$ on $(X,\hat g)$.
         Since if it were the case, according to Remark \ref{coneccc}
         the tangent cone of $S$ at $o$, a  {Simons} cone,
         would be calibrated in $(T_oX, \bar\Phi_o,\hat g_o)$.
         But $\bar\Phi_o$ can calibrate certain hyperplanes only. 
         Contradiction!
\\{\ }

Now we give a concrete construction for such $S$.
\\{\ }

\textbf{Example 2:}   
            Let $T$ be an oriented compact $(2r-1)$-dimensional smooth manifold.
            One can embed $K=S^{r-1}\times D^r$ into some small ball of $T$.
            After surgery along $S^{r-1}\times S^{r-1}$, one gets a manifold $T'$.
                   \begin{figure}[h]
\begin{center}
\includegraphics[scale=0.3]{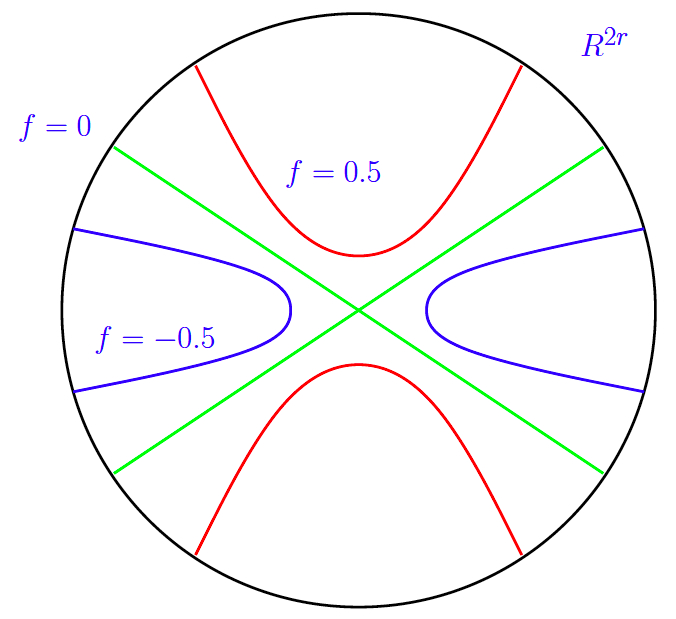}
\end{center}
\end{figure}
            The oriented $2r$-dimensional smooth manifold $W$
            obtained by the union of $[-0.5,0.5]\times (T-K)$ and the region between $\{f=-0.5\}$ and $\{f=0.5\}$ in the picture
            under the identification diffeomorphisms of $\{t\}\times S^{r-1}\times S^{r-1}$ with $\{f=t\}\bigcap S^{2r-1}(1)$
            provides a cobordism between $T$ and $T'$ (corresponding to $t=-0.5$ and $t=0.5$ respectively).
         Here $f$ is defined on $B^{2r}(1) \subset \mathbb R^r\times \mathbb R^r$   by 
         $f(\overrightarrow x, \overrightarrow y)=-\| \overrightarrow x\|^2+\|\overrightarrow y\|^2$,
         and $f^{-1}(0)$ is the truncated {Simons} cone.

          Take two copies of $W$. 
          Glue the same boundaries.
          Then one gets an orientable compact $2r$-dimensional manifold $X$.
          Now extend the Euclidean metric on the region between $\{f=-0.25\}$ and $\{f=0.25\}$ in the first copy to a metric on $X$.
          Let $(S,o)$ be the singular hypersurface in the first $W$ corresponding to $t=0$.
          Apparently $[S]\neq[0]$ in $H_{2r-1}(X;\mathbb R)$ (by intersection with a ``$t$-circle").
          Then Example 1 shows that 
          $S$ can be calibrated by a coflat calibration $\Phi$ singular only at $o$ with respect to some metric $g$ on $X$.
           
           \begin{rem}\label{88}
           By cross-products examples with more complicated singularity can be generated.          
           For instance, let $S_i,\Phi_i, X_i, g_i$ be given above for $i=1,2$.
           Then $S_1\times S_2$ with singularity $S_1\vee S_2$ is calibrated by the coflat calibration
           $\Phi_1\wedge \Phi_2$
           with singular set 
           $S_1\vee S_2$ 
           in the cartesian product $(X_1, g_1)\times (X_2, g_2)$.
           \end{rem}
            \begin{rem}
            Suppose $C$ is a $k$-dimensional cone in $ \mathbb R^{n}$ that has a calibration singular at most at one point.
            Consider $\Sigma_C\triangleq (C\times \mathbb R)\bigcap S^{n}(1)$ in $ \mathbb R^{n+1}$.
            Choose an $n$-dimensional oriented compact manifold $T$ with nontrivial $H_k(T;\mathbb R)$.
            Take an embedded oriented connected compact submanifold $M$ that represents a nonzero class of $H_k(T;\mathbb R)$.
            In smooth disks around a point of $M$ and a smooth point of $\Sigma_C$ respectively
            one can simultaneously connect $T$ and $S^n(1)$, $M$ and $\Sigma_C$
            through one surgery along $S^0\times S^n$ (i.e., connected sum).
            Denote by $X$ and $S$ the obtained manifold and submanifold (singular at two points).
            Then $[S]\neq 0\in H_k(X;\mathbb R)$ and similarly there exists a global calibration pair of $S$ by the proof of Theorem \ref{conecal}.
            \end{rem}

          \textbf{Example 3:} 
               Let $M$ be the smooth ``fiber" corresponding to  $\{t=-0.3\}$ in the first copy of $W$ in Example 2.
               Note that $\Phi$ is already a coflat calibration of $S$ on $(X, g)$.
               According to \S\ref{forms}, \S\ref{metrics} and Remark \ref{cptcst}, one can modify the calibration to $\tilde \Phi$
               and conformally change $g$ to $\tilde g$ in a neighborhood of $M$ away from $S$
such that $\tilde \Phi$ becomes a coflat calibration calibrating both $S$ and $M$ in $(X, \tilde g)$. 

However the homologically mass-minimizing submanifold $M$ cannot be calibrated by any smooth calibration in $(X, \tilde g)$. If it were,  then $S$ must be calibrated by the same smooth calibration as well which would lead to a contradiction.
This implies that all coflat calibrations of $M$
in $(X,\tilde g)$ share at least a common singular point.
For such creatures of higher codimension,
one can consider 
$M\times \{ \text{a point}\}$
in the Riemannian product of $(X,\tilde g)$ and a compact oriented manifold.
See Remark \ref{88} for more complicated examples.
\\

Next we consider the non-orientable case.
\\{\ }

\textbf{Example 4:}
Based upon $C_{3,4}$ one can get
an eight-dimensional oriented compact connected submanifold
$S$ with one
singular point in 
some oriented manifold $X^9$
with $[S]\neq [0]\in H_8(X;\mathbb R)$
by the method of Example 2.
Now blow up $X$ at a regular point of $S$.
Call the resulting manifold and submanifold $\check X$ and $\check S$ respectively.

By the {Seifert-van\ Kampen} theorem $\pi_1(\check X)\cong\pi_1(X)*\pi_1(\mathbb RP^8).$
The isomorphism of $\pi_1(\mathbb RP^8)\cong \mathbb Z_2$
trivially extends to a homomorphism $:$
$\pi_1(\check X)\rightarrow \mathbb Z_2$,
which canonically determines a two-sheeted cover $\overline X$ of $\check X$.
Denote the lifting of $\check S$ by $\overline S$.
Note that $\overline X\cong X\# X$ and $\overline S\cong S\#S^{\text{opposite orientation}}$.
By {\em Mayer-Vietoris} sequences, one has
$$H_8(\overline X;\mathbb R)\cong H_8(X;\mathbb R)\oplus H_8(X;\mathbb R), \text{and}$$
$$[\overline S]= [(S,-S)]\neq[0]\ \text{in} \ H_8(\overline X;\mathbb R).$$
Now create a $\mathbb Z_2$-invariant metric $\bar g$ on $\overline X$ such that
the orientable $\overline S$ can be calibrated (by a twisted calibration in the sense of \cite{TM}).

$\check S$ induces a $d$-closed integral current mod $2$,
$[[\check S]]_2$ (see \cite{Ziemer}), representing a nonzero $\mathbb Z_2$-homology class $[\check S]_2$.
We want to show that $[[\check S]]_2$ is $\mathrm{\mathbf M}^2$-minimizing in $[\check S]_2$ under the induced metric
$\check g$ on $\check X$,
where the \textbf{mass} $\mathrm{\mathbf{M}}^2(\cdot)$ of an integral current mod $2$
is the infimum of the mass of all integral representatives.

Suppose $K-[[\check S]]_2=dW$ in the sense of mod 2 for an integral current $K$ of finite mass and $W$ a top dimensional integral current mod 2.
Then the lifting expression to $\overline X$ becomes $\overline K-[[\overline S]]=d\overline W$ in the sense of mod $2$.
(Since $\overline S$ is orientable, $[[\overline S]]$ is an integral current up to a choice of orientation.)
Now $\overline X$ is oriented and $\overline W$ is of top dimension,
so $\overline W$ comes from
the quotient of $\tilde W$ by $2$
where $\tilde W$ is the integral current with multiplicity one on $\mathrm{\bold{spt}}(\overline W)$
and orientation inherited from $\overline X$.
Restrict $\tilde W$ to the connected component of $\bold{spt}(\overline W)$ to $\overline S$
and
denote it by
$\tilde W^\circ$.
Assign $[[\overline S]]$ the orientation induced from $\tilde W^\circ$.  
Let $-\overline K^\circ\triangleq d\tilde W^\circ-[[\overline S]]$.
It follows
\[
\mathrm{\mathbf M}_{\check g}(\check S)=
\frac{1}{2}\mathrm{\mathbf{M}}_{\bar g}([[\overline S]])\leq\frac{1}{2}\mathrm{\mathbf{M}}_{\bar g}(\overline K^\circ)\leq \mathrm{\mathbf{M}}_{\check g}(K).
\]
Running $K$ through all the integral representatives of $[[\check S]]_2$ one has
\[
\mathrm{\mathbf M}_{\check g}(\check S)=\mathrm{\mathbf M}_{\check g}^2([[\check S]]_2).
\]
Let $K_2$ be the integral current mod $2$ of an integral current $K$ with $[K_2]=[\check S]_2$.
Then
\[
\mathrm{\mathbf M}_{\check g}^2([[\check S]]_2)\leq \mathrm{\mathbf M}^2_{\check g}(K_2).
\]
Namely, $[[\check S]]_2$ is $\mathrm{\mathbf M}_{\check g}^2$-minimizing 
in its $\mathbb Z_2$ homology class.

\appendix

{\ }
\\
{\  }

\end{document}